\documentclass[11pt,a4paper,reqno]{article}
\usepackage{amsthm, mathrsfs,amssymb,amsmath,}
\usepackage{authblk}
\usepackage{epsfig}
\usepackage{graphicx,tikz}
\usetikzlibrary{arrows,automata}
\newtheorem{theorem}{Theorem}[section]
\newtheorem{lemma}[theorem]{Lemma}
\newtheorem{corollary}{Corollary}[theorem]
\theoremstyle{definition}
\newtheorem{definition}[theorem]{Definition}
\newtheorem{example}[theorem]{Example}
\newtheorem{proposition}[theorem]{Proposition}

\newtheorem{note}[theorem]{Note}
\theoremstyle{remark}
\newtheorem{remark}[theorem]{Remark}
\usepackage{hyperref}
\allowdisplaybreaks
\numberwithin{equation}{section}

\title{Fractal dimension for Inhomogeneous graph-directed attractors}

\usepackage{authblk}

\author[1]{Shivam Dubey}
\author[2*]{Saurabh Verma}
\affil[1,2*]{Department of Applied Sciences, IIIT Allahabad, Prayagraj 211015, India}	
\affil[1]{sd2461997@gmail.com}
\affil[2*]{correspondence to: saurabhverma@iiita.ac.in}

\begin{document}


\date{}
\maketitle





 

\begin{abstract}
In this paper, we define inhomogeneous Graph-Directed (GD) separation conditions for a given inhomogeneous GD Iterated Function Systems (IFS), and estimate the upper box dimension of attractors by the dimension of the condensation set and associated Mauldin-Williams graph dimension. Following some work of Fraser, we also estimate the lower box dimension of attractors generated by inhomogeneous GDIFS. In the end, we shed few lights on the continuity of dimensions for the attractors of inhomogeneous GDIFS.
\end{abstract}
{\bf Keywords:} {Graph-directed Iterated function systems, Inhomogeneous  GD attractors, Box dimension, Hausdorff dimension, Open set condition, Perron-Frobenious operators}\\
{\bf Math subject Classifaction:} {Primary 28A80; Secondary  28A70, 28A78, 47B65 }



\section{Introduction}\label{section 1} 
Mandelbort \cite{Mb} was the one who initially coined the term ``{\it Fractal"}. There is no univocal definition for the term `` {\it Fractal"}. In the literature, there are a bunch of tools available to construct a fractal set. Among all these tools, one of the methods to construct fractal sets is the ``{\it Iterated Function Systems "} (IFS) \cite{Hut}, which attracted many researchers for applications and constructions of fractals \cite{MF2}. Recently, many variants of IFS such as inhomogeneous IFS \cite{MF3}, Graph-Directed (GD) IFS \cite{Rd}, generalized GDIFS \cite{Nussbaum1}, and inhomogeneous GDIFS \cite{SS6} are introduced and studied in this field. The term inhomogeneous IFS introduced by Barnsley and Demko \cite{MF3},  which consists natural IFS with a compact set, called condensation set. Barnsley \cite{MF1} gives the idea of fractal interpolation function which is receiving a great attention, see, for instance, \cite{PAS3,EPJS,Chana,SS5,M2,SP}. Being an important index, fractal dimension is a vital part of fractal geometry \cite{Fal}. Dimension in the fractal geometry is always a key attraction to scientists and mathematicians. fractal sets can be defined, among other things, by its dimension. There are various fractal dimensions such as Minkowski-Bouligand dimension (box counting dimension) \cite{Fal}, Hausdorff dimension \cite{Fal}, Assouad dimension \cite{Fraserbook}, Packing dimension \cite{Fal}, quantization dimension \cite{GL1}, Higuchi dimension \cite{Higu}, entropy dimension \cite{MHOCHMAN}, $L^q$-dimension \cite{Shmerkin}, and many more. For studying and analyzing complex systems in various disciplines, such as mathematics, physics, biology, and computer science, a fractal dimension is a potent tool. It has various applications in different fields including modeling natural events, data analysis, and image compression. Often, exact dimension calculation in fractal geometry is an extremely challenging task, so many mathematicians estimate the dimension in terms of bounds, for example; Priyadarshi in \cite{Pr1} estimated a better lower bound on the Hausdorff dimension of a set of complex continued fractions. There are various techniques used in the estimation and calculation of fractal dimensions, such as $\delta$-covering technique \cite{Chana, Liang1}, operator theoretic technique (Perron-Frobenious operators) \cite{ Nussbaum1,Pr3, Pr1}, function space technique \cite{SS5, Chana, Verma21}, Fourier transformation technique \cite{Fal}, Potential theoretic technique \cite{Fal},  thermodynamical formalism \cite{RD1}, mass distribution principle \cite{Fal}, entropy \cite{MHOCHMAN}, and $L^q$-norms of convolution \cite{Shmerkin}, etc. Many researchers \cite{PAS3,EPJS,Chana,SS5,SP} studied the dimension of fractal functions using various methods. Using covering method, Liang and his collaborators \cite{Chana,Liang1} computed the exact value of the dimension of fractional integrals. Jha and Verma in \cite{Verma21} gave some bounds on the dimension of $\alpha$-fractal functions using different methods mentioned above. Recently, Chandra and Abbas \cite{SS5} try to give some estimates for the dimension of fractal functions by function space technique. Mauldin and Williams \cite{Rd} have given a more general GD construction. Fraser \cite{Fraser1} investigated the dimension of inhomogeneous attractors. Wang and Huang in \cite{WH1} studied the box dimension of inhomogeneous self-conformal sets with overlaps. In \cite{bFalc}, Boore and Falconer investigated the GD attractors and gave examples of attractors which cannot be generated by traditional IFS. Nussbaum, Priyadarshi, and Lunel in \cite{Nussbaum1} generalized the GD construction of \cite{Rd} and gave a formula for the Hausdorff dimension of invariant sets using the Perron-Frobenious operators. In \cite{Fraser1}, Fraser estimated that the countably stable dimension is the maximum of dimensions of the homogeneous attractor and the condensation set, also extended the results of \cite{Nina} by computing the upper box dimension, which behaves similar to countably stable dimensions. Baker, Fraser, and M\'ath\'e in \cite{Baker} obtained the new upper bounds for the upper box dimension of an inhomogeneous self-similar set with overlaps. Edgar and Golds \cite{EJG} estimated the upper bounds for the dimension of GDIFS attractors in terms of the dimension associated with Mauldin-Williams graph \cite{Rd}. Further, they also gave a lower bound on the dimension if each participating map in the IFS is bi-Lipschitz and the IFS satisfies the strong open set condition. 
\par

In this paper, we provide an upper bound for the dimension of the sets $F_i$ in Definition \ref{def1}. We use $\delta$-covering and operator theoretic techniques to estimate the fractal dimension of the inhomogeneous GD attractor. We should admit that our present work is partly influenced by Fraser \cite{Fraser1}, Priyadarshi \cite{Pr3}, Olsen and  Snigireva \cite{Nina}. We believe that the well-illustrated examples and remarks written in this paper will certainly benefit the non-expert and new readers to understand the dimension theory.
\par
The rest of this paper is assembled as follows. Section $2$ reviews the GDIFS. Section $3$ deals with the inhomogeneous GDIFS and the relationship between homogeneous GD fractal sets and inhomogeneous GD fractal sets. Section $4$ contains our main results. 

\section{A revisit to graph-directed systems and inhomogeneous Iterated Function Systems}\label{Section 2}

In the current section, we aim to give some previous definitions and results on graph-directed systems and inhomogeneous IFS. 

\begin{definition}
    Let $(X,d)$ be a complete metric space and $\mathcal{I} =\{X; f_1, f_2,\dots, f_N\}$ be an IFS such that $f_1,f_2, \ldots, f_N$ are contraction mappings on $X$. Then unique non-empty compact sets $A$ and $A_C$ are said to be the homogeneous attractor and inhomogeneous attractor with condensation set $C$ of the IFS $\mathcal{I}$, respectively, if $A$ and $A_C$ satisfy the following equations called  self-similar equations:
    $$A= \bigcup_{i=1}^N f_i (A) \text{ and } A_C = \bigcup_{i=1}^N f_i (A) \cup C, \text{ respectively}.$$ 
    Let $(p_1,\ldots,p_n)$ and $(p_1,\ldots,p_n,p)$ be probability vectors. Then Borel probability measures $\mu$ and $\mu'$ on \(X\) are said to be a homogeneous invariant measure and an inhomogeneous invariant measure with condensation measure $\nu$ with compact support $C$, respectively if they satisfy
    $$\mu= \sum_{i=1}^N p_i \mu \circ f_i^{-1}~\text{ and }~ \mu'= \sum_{i=1}^N p_i \mu' \circ f_i^{-1} + p \nu, \text{ respectively.}$$
\end{definition}

\begin{definition}
    If $S \subseteq X$ is a non-empty bounded subset of a complete metric space $(X,d)$. For $\delta >0$, let $N_{\delta}(S)$ be the minimum number of sets of diameter $\delta$ required to cover $S$ then upper box dimension of $S$ is defined as
    $$\overline{\dim}_{B}(S) = \limsup_{\delta \to 0^+ } \frac{\log N_{\delta}(S)}{-\log \delta},$$
    and lower box dimension of $S$ is defined as
     $$\underline{\dim}_{B}(S) = \liminf_{\delta \to 0^+ } \frac{\log N_{\delta}(S)}{-\log \delta},$$
     further, if both limits exist and are equal then the box dimension of $S$ is denoted by $\dim_{B}(S) $ and defined as
     $${\dim}_{B}(S) = \lim_{\delta \to 0^+ } \frac{\log N_{\delta}(S)}{-\log \delta}.$$
\end{definition}
\begin{definition}
   Let $S \subseteq X$ be a non-empty subset of a complete metric space $(X,d)$ and $r \in [0,\infty )$,
   $$H^{r}_{\delta}(S)=\inf \bigg\{ \sum_{i=1}^{\infty} |S_i|^r : S \subseteq \bigcup_{i=1}^{\infty} S_i, |S_i| < \delta \bigg\},$$
   where $|S_i|$ is the diameter of $S_i$ and the above infimum is taken over all countable cover of $S$. The $r$-dimensional Hausdorff measure of $S$ is defined as $H^r (S) = \lim_{\delta \to 0^+} H^{r}_{\delta} (S),$ and Hausdorff dimension of $S$ is defined by 
   $$\dim_H (S) = \inf\{ r \ge 0 : H^r (S) =0 \}.$$
 The Hausdorff dimension of a Borel probability measure $\mu$ is defined by 
 $$\dim_H(\mu)= \inf\{\dim_H(F): F ~\text{is a Borel subset such that } \mu(F)>0\}.$$
\end{definition}

\begin{definition}
 A directed graph $(V, E^*, i, t)$, where \(V\) is the set of all vertices and \(E^*\) is the collection of all directed paths, together with the initial vertex functions $i: E^* \rightarrow V $ and terminal vertex function $t: E^* \rightarrow V .$ 
 Let $e \in E^*$ be any finite path, then we may write $e= e_1e_2....e_k$ for some edges $e_i$. The initial vertex of $e$ is the initial vertex of its first edge $e_1$, so $i(e) = i(e_1)$  and similarly $t(e) = t(e_k).$ where $e_i \in E^1, ~~1\le i \le k$. $E^1$ is the set of all edges joining vertices $u~ \text{and}~ v $ for some adjacent vertices, $u ,v \in V.$ An edge $e \in E^*$ is said to be of length $n$ if $e = e_1e_2....e_n,$ where $e_i \in E^1.$
 \end{definition}

 \begin{remark}
For convenience, throughout this paper we will take the set of vertices of a directed graph as $V=(v_1,v_2....,v_N)$, also we will consider directed graph is strongly connected, for any two vertices \(u\) and \(v\) there is some directed edge joining them. 
\end{remark}
 \begin{definition}
 A GDIFS $(V, E^*,i,t,(X_{v_j},d_{v_j})_{v_j \in V},(f_e)_{e \in E^1})$ consists of $(V,E^*,i,t)$ a directed graph, and $(X_{v_j}, d_{v_j})_{{v_j} \in V}$ are complete spaces (each $X_{v_j}$ is a copy of $\mathbb{R}^n$ associated with vertex \(v_j\) ) and $f_e$ is contracting mapping from $X_{v_k}$ to $X_{v_j}$ with contracting ratio $\rho_e \in [0,1),$ such that $i(e)=v_j $ and $t(e)=v_k.$
 \end{definition}
For the construction of homogeneous GDIFS, one can visit \cite{Rd}, Boore, and Falconer \cite{bFalc} show the existence of attractors of GDIFS, which cannot be generated by standard IFS.
 
\section{Inhomogeneous GD fractal sets and measures}
We first revisit the inhomogeneous GD fractal sets and measures introduced in our recent work \cite{SS6}, and one may see that these sets and measures are natural generalizations of inhomogeneous fractal sets and measures appeared in \cite{MF3,Nina}. 

\begin{definition}
\label{def1}
Let $(V, E^*,i,t,(X_{v_j},d_{v_j})_{{v_j} \in V},(f_{e})_{e \in E^1})$ be a GDIFS and $(C_1,C_2,\ldots ,C_N) \in \prod_{j=1}^{N} K(X_{v_j})$ be fixed element, where $K(X_{v_j})$ is the space of all non-empty compact subsets of $X_{v_j}.$ A non-empty compact set  $(F_1, \ldots,F_N)$ such that

\begin{equation}
\label{eqn10}
(F_1,\ldots,F_N) =\bigg( \bigcup_{e \in E_1^{1}} f_{e}(F_{t(e)}) \cup C_1,\bigcup_{e \in E_{2}^1 } f_{e}(F_{t(e)}) \cup
C_2,\dots,\bigcup_{e\in E_{N}^1} f_{e}(F_{t(e)}) \cup C_N \bigg)
\end{equation}

is called an inhomogeneous GD Fractal set associated with the list $(f_e: e \in E_1, (C_1, C_2,..., C_N))$
\begin{note}
    In each union in (\ref{eqn10}) \(E_j^{1}\) is the set of all directed path \(e\) of length \(1\) with  \(i(e)= v_j\), and if \(e\) is a directed path joining vertex \(v_i\) to \(v_j\) then \(f_e\) is a map from $X_{v_j} \to X_{v_i}.$ 
    \end{note}
    \begin{note}
        If GDIFS has only one vertex then this construction reduces to the standard IFS and inhomogeneous GDIFS reduces to inhomogeneous IFS.
    \end{note}
\end{definition}
\begin{definition}\label{new612}
Let $(V, E^*,i,t,(X_{v_j},d_{v_j})_{{v_j} \in V},(f_{e})_{e \in E^1})$ be a GDIFS. Also let $\{(p^j_e)_{e \in E_j^1},p_j\}$ where $j \in (1,2,\ldots,N)$ be probability vectors corresponding to each $j$ and $\lambda_j$ be probability measures of compact support \(C_j\) on each $X_{v_j}$ associated to a vertex $v_j$, A probability vector-measure $(\mu_1,\ldots,\mu_N)$ such that
\begin{equation} \label{eqn12}
(\mu_1,\ldots,\mu_N) =\bigg (\sum_{e \in E_1^1} p^1_e \mu_{t(e)}\circ {f_e}^{-1} + p_1 \lambda_1,\ldots,\sum_{e \in E_N^1} p^N_e \mu_{t(e)}\circ {f_e}^{-1} + p_N \lambda_N \bigg)
\end{equation}
is called inhomogeneous GD invariant fractal vector-measure associated with the list $\{f_e: e \in E^1, ((p^j_e)_{e \in E_j^1},p_j),\lambda_j):j \in (1,2,\ldots, N)\}$
\end{definition}

\begin{proposition}\label{new614}
Let $(\mu_1, \ldots, \mu_N)$ be the fractal vector-measure  satisfying Equation (\ref{eqn12}) and let $C_i$ be the support of $\lambda_i$ for each $1 \le i \le N$. Then the support of $(\mu_1, \ldots, \mu_N)$ is a non-empty compact set $(F_1, \ldots, F_N)$ satisfying Equation (\ref{eqn10}).
\end{proposition}

\begin{theorem}
Let, $(O_1, \ldots,O_N)=\bigg(  \bigcup_{e \in E_1^*} f_e(C_{t(e)}) \cup C_1, \ldots, \bigcup_{e \in E_N^*} f_e(C_{t(e)})\cup C_N\bigg).$ Then
\begin{enumerate}
  \item[(i)] $(O_1, \ldots,O_N)$ satisfies inhomogeneous GD  equation(\ref{eqn10}).
 \item[(ii)] $(\overline{O}_1, \ldots,\overline{O}_N)=(F_1, \ldots, F_N)=(F_{\phi_1}\cup O_1, \ldots, F_{\phi_N} \cup O_N),$ where \((F_{\phi_1},\ldots, F_{\phi_N})\) is the corresponding homogeneous GD fractal set associated with the list $(f_e: e \in E_1, (C_1, C_2,\ldots, C_N)),$
 where $(O_1,\ldots,O_N) \in \prod_{j=1}^{N} K(X_{v_j}),$ such that
$$O_i= \bigcup_{e \in E_i^*} f_e(C_{t(e)}) \cup C_i, \hspace*{5em} \forall~ 1\le i \le N.$$ is orbital set.
\end{enumerate}
\end{theorem}

\section{Main Results}
In this section, we define GD inhomogeneous open set condition and strong open set condition for inhomogeneous GD system and also find a suitable upper bound for the upper box dimension of the attractor of an inhomogeneous GD system in terms of Mauldin-Williams graph dimension \cite{Rd}.

\begin{definition}
Let $\mathcal{GD}=(V,E^*,i,t,(X_{v_j},d_{v_j},C_j)_{v_j \in V},(f_e)_{e \in E^1})$ be a GDIFS, and 
$$(F_1, \ldots, F_N)=\bigg(  \bigcup_{e \in E_1^*} f_e(C_{t(e)}) \cup C_1, \ldots, \bigcup_{e \in E_N^*} f_e(C_{t(e)})\cup C_N\bigg),$$ is corresponding inhomogeneous GD  attractor, then it is said to be satisfying GD inhomogeneous open set condition (GDIOSC) if there exists a list of open sets $(U_1,\ldots, U_N) \in \prod_{j=1}^{N} K(X_{v_j}),$ such that
\begin{enumerate}
\item[(i)] $f_e(U_{t(e)}) \subset U_j$ for each $e \in E_j ^1~ \text{and} ~ j \in \{1,2 \ldots, N \}$
\item[(ii)] For each $ e \neq e',$
$f_e(U_{t(e)}) \cap f_{e'}(U_{t(e')}) = \phi.$
\item[(iii)] For each $j \in (1, \ldots, N),$
$C_j \subseteq \overline{U}_j.$
\item[(iv)] If  for each $j \in (1, \ldots, N),$ ~~ $U_j \cap F_j \neq \phi$ then it is said to be satisfying GD Inhomogeneous strong open set condition (GDISOSC).

\end{enumerate}  
\end{definition}
\begin{remark}
Here we give some examples to understand the inhomogeneous open set condition (OSC).
\begin{itemize}
    \item Let an IFS containing two functions $S_1(x)= \frac{x}{3}$ and $S_2(x)= \frac{x}{3} + \frac{2}{3}$ on $\mathbb{R}.$ This IFS satisfies the OSC with open set $(0,1),$ also the corresponding inhomogeneous system with condensation set $C = \bigg[ \frac{1}{3} , \frac{2}{3}\bigg]$ satisfies the inhomogeneous OSC with same open set $(0,1).$
    \item Let $S_1, S_2: \mathbb{R} \to \mathbb{R}$ be defined by $S_1(x)= \frac{x}{2} ~~\text{and}~~ S_2= \frac{x}{2}+ \frac{1}{2},$ then clearly the IFS with these two functions satisfies the OSC with open set $(0,1).$ But the corresponding inhomogeneous IFS with condensation set $C= \{2\}$ does not satisfy the inhomogeneous OSC with the same open set $(0,1).$
    \item Let $\{S_1,S_2 : \mathbb{R}^2 \to \mathbb{R}^2\}$ be an IFS with $S_1(x,y)= \frac{1}{\sqrt{2}} (x,y) ~~ \text{and} ~~ S_1(x,y)= \frac{1}{\sqrt{2}} (x,y) + (1-\frac{1}{\sqrt{2}}, 0).$
The homogeneous attractor associated with the IFS is $F_{\phi} = [0,1] \times \{0\}$ with $\overline{\dim}_{B} (F) = 1,$ taking condensation set $C =\{0\} \times [0,1]$ with $\overline{\dim}_{B} (C) = 1,$ then corresponding inhomogeneous IFS $F_{C}$ has dimension $\overline{\dim}_{B} (F_{C}) = \frac{\log 8}{\log 4} > 1$ , which does not follow $\overline{\dim}_{B} (F_{C}) \le \max \{ \overline{\dim}_{B} (F_{\phi}), \overline{\dim}_{B} (C)\}.$ This IFS does not satisfy inhomogeneous OSC, see \cite{Baker} for more details.
\end{itemize}
\end{remark}

\begin{definition} 
Let $\mathcal{GD}=(V,E^*,i,t,(X_{v_j},d_{v_j},C_j)_{v_j \in V},(f_e)_{e \in E^1})$ be a GDIFS. Let $\{\mathcal{M}(s)\}_{s \ge 0}$ be $N \times N$ matrices indexed by vertices of the graph associated with GDIFS, with entries
$\mathcal{M}_{v_i v_j}(s)=\sum_{e \in E_{ij}^1} \rho_e^s.$
Let $\Phi$ be a function defined as 
$\Phi (s)= \Upsilon(\mathcal{M}(s)).$  Where $\Upsilon(\mathcal{M}(s))$ is spectral radius of $\mathcal{M}(s).$ It is a strictly decreasing, continuous function with $\Phi(0) \ge 1$, also $\lim_{s \to \infty }\Phi(s) = 0.$ A unique positive real number $s^*$ such that $\Phi(s^*) = 1$ is the Mauldin-Williams graph dimension \cite{EJG} associated with the GDIFS $\mathcal{GD}.$
\end{definition}
Clearly, $\{\mathcal{M}(s^*)\}$ is an irreducible matrix with positive entries and $1$ is an eigenvalue so by Perron-Frobenius theorem there must be a positive eigenvector $(u_1,u_2,\cdots,u_N)$:
\begin{equation}\label{eq90}
  u_i > 0,~~~ \sum_{i=1}^{N} u_i = 1, ~~~~~\sum_{1 \le j \le N} \sum_{e \in E_{ij}^1} \rho(e)^{s^*} u_j = 1.u_i, ~~~~~~~ 1\le i \le N  
\end{equation}
\begin{example}
    Here we consider an irreducible matrix with positive entries 
    $$\begin{pmatrix}
  \frac{1}{2} & \frac{1}{3} \\ 
  ~~\\
  \frac{1}{2} & \frac{1}{3}
\end{pmatrix}.$$
This matrix has eigenvalues $0 ~~\text{and} ~~\frac{5}{6}.$ Clearly its spectral radius $\frac{5}{6} $ is one of its eigenvalues and there is a positive eigenvector $\Big(\frac{1}{2} ~~~ \frac{1}{2}\Big)$ corresponding to eigenvalue $\frac{5}{6}$ with $\frac{1}{2} + \frac{1}{2} =1.$
\end{example}


\begin{example}
    Let us consider a GD system consisting of two vertices $v_1$ and $v_2$ and five directed edges, namely $e_1,e_2,e_3,e_4$ and, $e_5$ as shown in Figure 1. Let vertices $v_1$ and $v_2$ associated with usual metric space $X_{v_1}=X_{v_2}=(\mathbb{R},d)$ and directed edges associated with the contraction mappings given below
$f_{e_1}= \frac{x}{2}, ~ f_{e_2}= \frac{x}{4},~ f_{e_3}= \frac{x}{6},~f_{e_4}= \frac{x}{8},~f_{e_5}= \frac{x}{10}.$
For any $s>0$, the matrix associated with the GD system is 
$\mathcal{M}(s)=\begin{bmatrix}
\frac{1}{4^s} & \frac{1}{2^s}+\frac{1}{6^s} \\ 
\frac{1}{10^s} &  \frac{1}{8^s}
\end{bmatrix}$
\begin{figure}
\centering
\begin{tikzpicture}[->,>=stealth',shorten >=1pt,auto,node distance=5cm,
                    thick, main node/.style={circle,draw,font=\sffamily\small\bfseries}]

 \node[main node] (v1) {v1};
 
  \node[main node] (v2) [right of= v1] {v2};

  \path[every node/.style={font=\sffamily\large}]
   
    (v1) edge node [right] {$e_1$} (v2)
      
        edge [loop left] node {$e_2$} (v1)
        edge [bend right] node[left] {$e_3$} (v2)
    
    (v2) edge [loop right] node {$e_4$} (v2)
        edge [bend right] node[right] {$e_5$} (v1);
\end{tikzpicture}
 \caption{Graph-directed system with $2$ vertices and $5$ edges}
\end{figure}
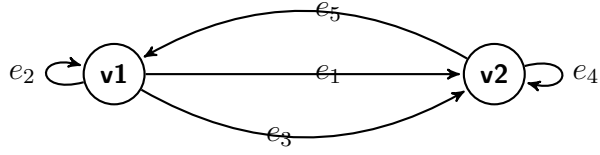
\end{example}
\begin{remark}\label{new189}
   Through this remark, we aim to use the method of positive operators (widely known as the Perron-Frobineous operators) due to Priyadarshi and his group \cite{Nussbaum1} for the estimation of dimension. First define $$\mathcal{C}(X_{v_j}):=\{g: X_{v_j} \to \mathbb{R} : g~ \text{is continuous and bounded}\},$$ for $1 \le i \le N$ with sup norm, $\|g\| = \sup_{x \in X_{v_j}} |g(x)|.$

  For $s>0,$ consider a linear map $F_s : \prod_{i=1}^N \mathcal{C}(X_{v_i}) \to \prod_{i=1}^N\mathcal{C}(X_{v_i})$ defined by    \[F_s(g_1,\cdots,g_N)(x_1, \cdots, x_N) = \bigg( \sum_{ e \in E,~ t(e)=v_1} \rho_e^s g_{i(e)}(f_{e}(x_1)), \ldots, \sum_{ e \in E,~ t(e)=v_N} \rho_e^s g_{i(e)}(f_{e}(x_N)) \bigg),\]
   where  $x_i \in X_{v_i}.$
\end{remark}
\begin{remark}
With the same notation as in \cite{Nussbaum1}, if we consider $\Gamma(\alpha, e) = i(e)$ then the generalized GD construction of \cite{Nussbaum1} will reduce to the Mauldin-Williams GD system, see, for instance, \cite[Example 2.1]{Nussbaum1}.
\end{remark}
\begin{remark}\label{new143}
   By \cite[Lemma 4.6]{Nussbaum1}, there exists a unique $s >0$ with spectral radius of $F_s$ is $1,$ and a non-zero eigenvector $g=(g_1, \ldots,g_N)$ such that $\|g\|=1.$
\end{remark}

\begin{definition}[\cite{EJG}]
    Let $E^{(\infty)} \subseteq E^*$ be the collection of all paths having length infinite, then a cross-cut $T$ is the collection of all paths $e=e_1e_2 \ldots$ for which there exists only one $n$ such that $e_{|_n}= e_1\ldots e_n \in T.$ also for $1\le i,j \le N,~~~~ T_i= T \cap E_i^*$ and $T_{ij} = T \cap E_{ij}^*.$
\end{definition}
\begin{note}
    For $e = e_1, \ldots, e_n, \ldots$ we define $e^- $ by removing the last edge in $e.$
\end{note}
\begin{lemma}

 Let $T_i$ be $i^{th}$ cross-cut, then for any $\delta >0,$ ~~
 \begin{equation*}
\label{eq1}
|T_i| \le (\delta \rho_{min})^{-s^*} \frac{u_{max}}{u_{min}},
\end{equation*}
where $\rho_{min} = \min_{e \in E_i}\{\rho(e)\}$, $u_{\max} = \max_i u_i ,$ and $u_{\min}=\min_i u_i. $
\end{lemma}
\begin{proof}
    For any fixed $\delta >0$, we have $i^{th}$ cross-cut is given by
    $$T_i = \{ e : ~~~e \in E_{ij}^*,~~~~ \rho(e) < \delta \le \rho(e^-)\}.$$
    If $\rho_{min} = \min_{e \in E_i}\{\rho(e)\}$  then for any $e \in T_i,$ we have $\delta \rho_{\min} \le \rho(e) < \delta.$ Using the fact that $T_i$ constitutes a partition of $E_i^{\infty}$ and $\mu_i$ defined as $\mu_i([e])= \rho(e)^{s^*} u_j $  extended to a Borel measure on $E_i^{*}$ (for more details, see \cite{Rd}) then by equation (\ref{eq90}) we have 
   \begin{equation}\label{eqn616}
       \begin{aligned}
            u_i &= \sum_{j \in V} \sum_{e \in T_{ij}} \rho(e)^{s^*} u_j \\& \ge (\delta \rho_{\min})^{s^*} u_{\min} |T_i|.
       \end{aligned}
   \end{equation}
   Since $u_i>0$ (see equation (\ref{eq90})), it follows that $u_{min}>0$. This with equation (\ref{eqn616}) yields
   $$|T_i| \le (\delta \rho_{min})^{-s^*} \frac{u_{max}}{u_{min}},$$ which is the required result.
    
\end{proof}
\begin{theorem}[\cite{EJG}]
	Let $s^*$ be the Mauldin-Williams graph dimension associated with the GD system $\mathcal{GD}$, and let $(F_{\phi_1},\ldots,F_{\phi_N})$ be homogeneous GD fractal set, then $\dim_{H}F_{\phi_i} \le s^*.$ 
\end{theorem}
 \subsection{Upper box dimension} 
 In this subsection, we first obtain an upper bound for the upper box dimension of the orbital sets $O_i$. Then using some properties of dimension, we will give a suitable bound for the upper box dimension of GD inhomogeneous attractors $F_i.$
 \begin{theorem}\label{new613}
If $\overline{\dim}_{B}(C_i) = \overline{\dim}_B (C_j) ~ \text{for every} ~~ 1 \le i , j  \le N,$    then $
\overline{\dim}_B (O_i) \le \max\{s^*,\overline{\dim}_B (C_i)\}, ~~ \forall ~~1 \le i \le N.$     \end{theorem}
 \begin{proof}
 Let us begin by considering  $t > \max \{s^*, \overline{\dim}_B (C_i) \}$, that is, $\lim_{\delta \to 0}\frac{\log N_{\delta} (C_i)}{-\log {\delta}}  < t.$ Then by the definition of upper box dimension, we can choose a constant $c_i^t> 0$ such that $N_{\delta} (C_i) \le c_i^t \delta^{-t},$ holds for all sufficiently small $ \delta>0.$
Since, this is true for every $1 \le i \le N$. We have
\begin{equation*}
    \begin{aligned}
        N_{\delta}(O_i) &=  N_{\delta}\big(C_i \cup \bigcup_{e \in E_i^*} f_e(C_{t(e)} )\big) \\& \le  \sum_{e \in E_i^*, ~ \delta \le \rho(e)} N_{\delta}(f_e(C_{t(e)}))+ N_{\delta} \big(\bigcup_{e \in E_i^*, ~ \delta > \rho(e)} f_e(C_{t(e)})\big) +N_{\delta} (C_i) \\& \le \sum_{e \in E_i^*, ~ \delta \le \rho(e)} N_{\delta / \rho(e)} (C_{t(e)}) + \sum_{e \in T_i} N_{\delta / \rho(e)}(C_{t(e)}) + c_{\max}^t \delta^{-t} \\&  \le \sum_{e \in E_i^*, ~ \delta \le \rho(e)} N_{\delta / \rho(e)} (C_{t(e)}) + \sum_{e \in T_i} N_1(C_i) + c_{\max}^t \delta^{-t} \\&  \le c_{\max}^t \delta^{-t} \sum_{e \in E_i^*, ~ \delta \le \rho(e)} \rho (e)^t +|T_i| N_1^{\max}(C_i) + c_{\max}^t \delta^{-t} \\& \le c_{\max}^t \delta^{-t} \sum_{e \in E_i^*, ~ \delta \le \rho(e)} \rho (e)^t + \delta^{-s_1} (\rho_{\min})^{-s^*} \frac{u_{\max}}{u_{\min}}N_1^{\max}(C_i) + c_{\max}^t \delta^{-t} 
        \\& \le c_{\max}^t \delta^{-t} \sum_{e \in E_i^*, ~ \delta \le \rho(e)} \rho (e)^t + \delta^{-t} (\rho_{\min})^{-s^*} \frac{u_{\max}}{u_{\min}}N_1^{\max}(C_i) + c_{\max}^t \delta^{-t} 
        \\& \le (c_{\max}^t p_t + (\rho_{\min})^{-s^*} \frac{u_{\max}}{u_{\min}}N_1^{\max}(C_i) + c_{\max}^t) \delta^{-t}.
    \end{aligned}
\end{equation*}
The second inequality in the above arrives by choosing $c_{\max}^t= \max_{i} \{c_i^t\}$, and the third inequality comes due to the compactness of each $C_i: 1 \le i \le N,$ the fourth inequality follows by choosing $N_1^{\max}(C_i) = \max_i\{N_1(C_{t(e)})\}$, as 
$C_i$'s are compact sets, so each $N_1(C_i)$ is a fixed finite quantity, and the fifth inequality is due to the above lemma. Consequently, $\overline{\dim}_B (O_i) \le t.$ Now, since $t$ is any number arbitrarily close to $\max\{s^*, \overline{\dim}_B (C_i)\}$, we have $\overline{\dim}_B (O_i) \le \max\{s^*,\overline{\dim}_B (C_i)\}.$ This completes the proof of the theorem.
 \end{proof}
\begin{corollary}
$\max\{\overline{\dim}_B (F_{\phi_i}),\overline{\dim}_B (C_i)\} \le  \overline{\dim}_B (F_i) \le \max\{s^*, \overline{\dim}_B (C_i)\},$ For all $~~1 \le i \le n,$ 
\end{corollary}
\begin{proof}
Since $F_i = F_{\phi_i} \cup O_i$ and the upper box dimension is finitely stable, we have
$$\overline{\dim}_B (F_i) = \max\{\overline{\dim}_B (F_{\phi_i}),\overline{\dim}_B (O_i)\}.$$
 Now by using the monotonicity property of the upper box dimension and the above result, we get
 $$\max\{\overline{\dim}_B (F_{\phi_i}),\overline{\dim}_B (C_i)\} \le  \overline{\dim}_B (F_i) \le \max\{s^*, \overline{\dim}_B (C_i)\},$$
 establishing the claim.
\end{proof}
\begin{corollary}\label{result20}
    If all the mappings $f_e$'s in GDIFS are similarity transformations, and it satisfies GDISOSC then $$ \overline{\dim}_B (F_i) = \max\{s^*, \overline{\dim}_B (C_i)\}.$$
\end{corollary}
\begin{proof}
   Clearly under the assumption of GDISOSC along with all $f_e$'s are similarity transformation 
   $\overline{\dim}_B (F_{\phi_i})= s^*,$ and by the above corollary we get the result. 
\end{proof}
\begin{proposition}
    Under the assumptions of Theorem \ref{new613} and Definition \ref{new612}, we have $\dim_H(\mu_i) \le \max\{s^*, \overline{\dim}_B (C_i)\},$ for all $1 \le i \le N.$
\end{proposition}
\begin{proof}
    Using Proposition \ref{new614}, Theorem \ref{new613} and definition of the Hausdorff dimension of a measure, the results follow.
\end{proof}
\begin{remark}
    Using \cite[Theorem 4.17]{Nussbaum1} and \cite[Theorem 4.7]{sverma4}, Remark \ref{new189} produces that the Mauldin-Williams graph dimension will be the same as the number obtained in Remark \ref{new143} with the assumption of GDISOSC.
\end{remark}
\begin{proposition}\label{new617}
   If all mappings $f_e$'s of GDIFS $\mathcal{GD}$ satisfy
   \begin{equation}\label{eq90}
     d(f_e(x),f_e(y)) \ge \rho_e' d(x,y),
   \end{equation}
 
   for all $x,y \in X_{v_i}, ~ e \in E^1 , t(e)=v_i$ and $\mathcal{GD}$ satisfies  the GDISOSC. Then $s' \le \dim_H (F_i) \le \underline{\dim}_B(F_i), $ where $s'$ is the corresponding Mauldin-Williams graph dimension.
\end{proposition}
\begin{proof}
   It is well-known \cite{Fal} that $\dim_H (F_i) \le \underline{\dim}_B(F_i)$. One may easily prove that the corresponding homogeneous GDIFS will satisfy GDSOSC because of the assumption of GDISOSC. Let us now recall \cite[Theorem 3.9]{Rd} that  $s' \le \dim_H (F_{\phi_i})$ for all $1 \le i \le N$ under the GDSOSC. Since $F_i=F_{\phi_i} \cup O_i$, the previous two statements with monotonicity of Hausdorff dimension yield $s' \le \dim_H (F_i)$, completing the assertion.
\end{proof}
\begin{remark}
 All similarity transformations satisfy the assumption of the above proposition. However, in this remark, we give some examples which are not similarity transformations.
First, consider $f:[\cos 1, 1] \to [\cos 1, 1]$ as $f(x)=\cos(x)$. By the mean-value theorem, we have $$\sin(\cos 1) |x-y| \le |f(x)-f(y)| \le|x-y|, \text{ for all } x, y \in [\cos 1, 1].$$ 
Next, for $0<m <n,$ define $f: \mathbb{R}^2 \to \mathbb{R}^2 $ as $$f(x,y)=\Big(\frac{x}{n}+\frac{1}{n},\frac{y}{m}+\frac{1}{n}\Big). $$ Note that 
 \[
 \frac{1}{n} \|(x,y)-(z,w)\|_2 \le \|f(x,y)-f(z,w)\|_2 \le \frac{1}{m} \|(x,y)-(z,w)\|_2,
 \]
 holds for all $(x,y), (z,w) \in \mathbb{R}^2.$
\end{remark}

\begin{theorem}\label{new615}
 If $\dim_{H}(C_i) = \dim_H (C_j)$ for every $1 \le i , j  \le N,$ and all $f_e$'s are similarity transformations then 
\begin{itemize}
    \item  $\dim_H (F_i) = \max \{\dim_H (F_{\phi_i}), \dim_H (C_i)\}$,
\end{itemize}
Further, if the GDIFS satisfies GDISOSC then
\begin{itemize}
    \item  $\dim_H (F_i) = \max \{\dim_H (F_{\phi_i}), \dim_H (C_i)\}=\max \{s^*, \dim_H (C_i)\}$,
\end{itemize}
\end{theorem}
\begin{proof}
    It is well-known \cite{Fal} that the Hausdorff dimension is countably stable and invariant under bi-Lipschitz map. This fact with the assumption $\dim_H (C_i) = \dim_H (C_j)$ for all $i$ and $j$ produces
    $ \dim_H (O_i)= \dim_H (C_i)$. Since $F_i= F_{\phi_i} \cup O_i$, it follows that $\dim_H (F_i) = \max \{\dim_H (F_{\phi_i}), \dim_H (C_i)\}$. Now, Proposition \ref{new617} gives the second part, this completes the proof.
    \end{proof}

\subsection{Lower box dimension}
In this subsection, we try to estimate some non-trivial bounds for the lower box dimension similar to \cite{Fraser1} under a certain assumption and using the concept of covering regularity exponent (CRE). Clearly, If all the mappings in GDIFS $\mathcal{GD}$ satisfying Equation $\ref{eq90}$, then $s'$ is a trivial lower bound for lower box dimension.

For a bounded set $S$ and for $0 <\delta \le 1, 0 \le t$, 
$$P_{t,\delta}(S)= \sup_{[0,1]}\{p: N_{\delta^p} \ge \delta^{-pt}\}$$ is $(t,\delta)$-CRE of $S,$ and the quantity $P_t(S) = {\lim \inf}_{\delta \to 0} P_{t, \delta} (S)$ is called $t$-CRE of $B.$

\begin{example}
If $B$ is a finite set:
For $t=0,~~ \delta^{-pt}=1$ ~~and ~~ $N_{\delta^p} \ge 1$ ~~ for all $0\le p\le 1.$ Also if $S$ is a finite set with $m$ elements, then there exists some $0 < \delta_0,$ such that for all $\delta \le \delta_0, ~~N_{\delta^p} = m$ for any $p \in (0,1].$

\[P_t(B)= \left\{\begin{array}{cl} 1, 
& t = 0 \\ 0,& t \ge 0 \end{array} \right.\]
If $S=[0,1],$ then
\[P_t(B)= \left\{\begin{array}{cl} 1, 
& t \le 1 \\ 0,& t >1 \end{array} \right.\]

If $t =0$ then clearly $P_t(S) =1$. If $ 0 \le t \le 1$ then the minimum number of boxes required to cover $S$ of diameter $\delta^p$ is $\big[ \frac{1}{\delta^p} \big]$, and for all $0< p \le 1$ we have 
$$\bigg[ \frac{1}{\delta^p} \bigg] \ge \frac{1}{\delta^{pt}},$$ and
$$\bigg[ \frac{1}{\delta} \bigg] \ge \frac{1}{\delta^{t}}.$$
also if $t > 1$ then obviously $P_t(S) =0.$
\end{example}
From the definition of $P_t (S)$ it is clear that it is a decreasing function in $t$ and moreover, we have the following lemma which is borrowed from \cite{Fraser1}. It is worth to note that the proof of this lemma is obvious by the definition of the upper box dimension and the lower box dimension of $S.$.
\begin{lemma}
  We have 
  \[ 
 P_t(B)=\left\{\begin{array}{cc}
    1 &  t > \overline{\dim}_B (S) \\
    0 &   t < \underline{\dim}_B (S)
\end{array}
\right.
.\] 
\end{lemma}
\begin{lemma}[\cite{Fraser1}]\label{eq2}
We have the following:
\begin{itemize}
    \item $0 \le P_{t,\delta}(C_i)$ and $P_t(C_i) \le 1,$ for all $t, \delta >0. $
    \item  $ N_{\delta^{P_{t,\delta}(C_i)}}(C_i) \ge \delta^{-P_{t,\delta}(C_i) t},$ for all $\delta>0.$
\end{itemize}
    \end{lemma}
    \begin{definition}
       Let $(X,d)$ be a metric space, and $\dim_H(X)= q < \infty $ then it is said to be $q-$Ahlfors regular if there exists a constant $c>0$ such that,
       $$c^{-1} r^ {\dim_H (X)} \le \mathcal{H}^{\dim_H (X)}(B(x,r)) \le c r^{\dim_H (X)},$$ for all, $x \in X~~~ \text{and}~~ 0 \le r \le \text{diam} (X). $
    \end{definition}
    \begin{example}
The set $B=[0,1],$ is $1-$Ahlfors regular. Self-similar subsets of $\mathbb{R}^n$
satisfying the open set condition
are standard examples of regular sets.
\end{example}
    \begin{lemma}[\cite{Fraser1}]
    Let $(X,d)$ be an Ahlfors regular space and let $(U_i)$ be a collection of disjoint open subsets such that each $U_i$ contains an open ball of radius $ar$ and contained in an open ball of radius $br$, then any open ball of radius $r$ can intersect $U_i$ at most 
    $c^2 \big(\frac{1+2b}{a}\big) ^ {dim_H (X)}.$
    \end{lemma}
    To increase the readership of this paper and for the benefit of new researchers, we include the following remark.
    \begin{remark}
       \begin{itemize}
           \item The set $\mathbb{Q} \cap [0,1]$ is not $0$-Ahlfors regular.
           \item A compact subset $X$ of $\mathbb{R}$ such that $\dim_H(X)=1$ and $\mathcal{H}^1(X)=0$ will not be Ahlfors regular space. Let us see how to construct such a compact set. For each $r<1$, one may construct a Cantor-type set with Hausdorff dimension $r$ by varying the lengths of the intervals in the usual Cantor set construction. In particular, we define $X_n \subset [0,1]$ as a Cantor-type set of Hausdorff dimension $1 - \frac{1}{n}$ for each $n \ge 2$. The union $X=\cup_{n =2}^{\infty} X_n$ then has Lebesgue measure $0$ because each $X_n$ does so, but Hausdorff dimension $1$. Since $\mathcal{H}^1(X)=0$, $\mathcal{H}^1(B(x,r))=0$ for all $x \in X$ and $r>0$. So, we cannot find any real number $c>0$ such that $$c^{-1} r^ {1} \le \mathcal{H}^{1}(B(x,r))=0.$$ Therefore, $X$ is not Ahlfors regular space.
           \item Note also that if $(X,d)$ is not Ahlfors regular space then the conclusion of the above lemma may not hold.
       \end{itemize}
    \end{remark}

\begin{theorem}
   If GDIFS be chosen such that $(X_{v_i},d_i)$ is Ahlfors regular and all $C_i=C_j$ and satisfies the GDIOSC, then
   $P_t(C_i) t + (1 - P_t(C_i))s' \le \underline{\dim}_B (F_i), $ for all $t \ge0.$
   
\end{theorem}
\begin{proof}
Define the $i^{th}$-crosscut $T_i$ for $\delta^{1-P^{t,\delta}(C_i)}{\rho'}_{\min}^{-1} >0$,

$$T_i=\{e: e \in E_{ij}^*, ~~\rho'_e < \delta^{1-P_{t,\delta}(C_i)}{\rho'}_{\min}^{-1} \le {\rho'}_{e^-}\}. $$ Now we can estimate the cardinality of the set $T_i$ as below:
\begin{equation*}
\begin{aligned}
  u_i &= \sum_{j \in V} \sum_{e \in T_{ij}} {\rho'}_e^{s'} u_j \\& \le N |T_i| (\delta^{1-P_{t,\delta}(C_i)}{\rho'}_{\min}^{-1} )^{s'}u_{\max}.
  \end{aligned}
\end{equation*}
This in turn yields
$$|T_i| \ge \frac{u_{\min}}{N u_{\max}} (\delta^{1-P_{t,\delta}(C_i)}{\rho'}_{\min}^{-1} )^{-s'}.$$
Let $t> \max\{s', \underline{\dim}_B (C_i)\}.$ Consider the case when $0 < \delta^{1-P^{t,\delta}(C_i)}{\rho'}_{\min}^{-1} \le 1$ and
let $(U_1,\cdots, U_N)$ be the list of open sets in the GDIOSC, then we can find two positive real numbers $m_i$ and $n_i$ for each $i$ such that $U_i$ contains an open neighborhood of radius $m_i$ and contained in an open neighborhood of radius $n_i.$ If $m$ is $\min(m_i)$ and $n$ is $\max(n_i)$ then for $e \in T_i,$ the set $f_e(U_i)$ is an open set which contains an open neighborhood of radius $m \delta^{1-P_{(t,\delta)}(C_i)}$ and is contained in an open neighborhood of radius $n {\rho'}_{\min}^{-1}$ also due to the GDIOSC, $\{f_e(U_{t(e)}): e \in T_i\}$ is family of pairwise disjoints sets in $X_{v_i}$. Since $f_e(C_i) \subseteq \overline{f_e(U_i)}$ and $X_{v_i}$ is Ahlfors regular, so any set of diameter $\delta >0$ cannot intersect more than $\lambda^2\Big( \frac{1+2n \rho_{\min}^{-1}}{m} \Big)^{\dim_H (X_{v_i})}$ of the sets $\{f_e(C_{t(e)}): e \in T_i\}.$
With these previous ingredients, we may proceed as follows:
\begin{equation*}
    \begin{aligned}
     N_\delta(O_i) &= N_\delta \big(C_i \cup \bigcup_{e \in E_i^*} f_e(C_{t(e)}) \big) \\& \ge \frac{\sum_{e \in T_i} N_{\delta}(f_e(C_{t(e)}))}{\lambda^2\big( \frac{1+2n \rho_{\min}^{-1}}{m} \big)^{\dim_H X_{v_i}}}  \\& \ge \frac{\sum_{e \in T_i} N_{\delta / \rho'_e}(C_{t(e)})}{\lambda^2\big( \frac{1+2n \rho_{\min}^{-1}}{m} \big)^{\dim_H X_{v_i}}} \\& \ge \frac{\sum_{e \in T_i} N_{\delta^{P_{t,\delta}(C_i)}}(C_{t(e)})}{\lambda^2\big( \frac{1+2n \rho_{\min}^{-1}}{m} \big)^{\dim_H X_{v_i}}} \\& \ge \frac{\delta^{-t P_{t,\delta}(C_i)}|T_i|}{\lambda^2\big( \frac{1+2n \rho_{\min}^{-1}}{m} \big)^{\dim_H X_{v_i}}} \\& \ge \frac{\delta^{-t P_{t,\delta}(C_i)} (\delta^{1-P^{t,\delta}(C_i)}{\rho'}_{\min}^{-1})^{-s'} u_{\min}/ (N u_{\max})}{\lambda^2\big( \frac{1+2n \rho_{\min}^{-1}}{m} \big)^{\dim_H X_{v_i}}} \\& = \frac{{\rho'^{s'}_{\min}}u_{\min}/(N u_{\max}) \delta^{-(P_{t,\delta}(C_i)t+(1-P_{t,\delta}(C_i))s')}}{\lambda^2\big( \frac{1+2n \rho_{\min}^{-1}}{m} \big)^{\dim_H X_{v_i}}} \\& \ge \frac{{\rho'^{s'}_{\min}}u_{\min}/(N  u_{\max} )\delta^{-(P_t(C_i) -\epsilon )t+(1-(P_{t}(C_i) - \epsilon))s'}}{\lambda^2\big( \frac{1+2n \rho_{\min}^{-1}}{m} \big)^{\dim_H X_{v_i}}}
    \end{aligned}
\end{equation*}
Third inequality comes as $\delta^{1-P^{t,\delta}(C_i)}{\rho'}_{\min}^{-1} \le 1$ this implies that $\delta {\rho'}_{\min}^{-1} \le \delta^{P^{t,\delta}(C_i)}$ ~~ and
the last inequality comes as for small $\epsilon \in (0,P_{t,\delta}(C_i)],$ we have $N_{\delta^{P_{t,\delta}(C_i)}}(C_i) \ge N_{\delta^{P_{t,\delta}(C_i)- \epsilon}}(C_i) \ge \delta^{-(P_{t,\delta}(C_i)- \epsilon)t}.$ 
This shows that $ \underline{\dim}_B (O_i) 
\ge (P_t (C_i) - \epsilon)t + (1- (P_t (C_i) - \epsilon))s'.$
Also when $\delta^{1-P^{t,\delta}(C_i)}{\rho'}_{\min}^{-1} \ge 1,$ we have $\delta^{-(1-P_{t,\delta}(C_i))^{s'}}{\rho'}_{\min}^{s'} \le 1.$ This implies that, 
$$N_{\delta}(O_i) \ge N_{\delta^{P_{t,\delta}(C_i)}}(C_i) \ge \delta^{-(1-P_{t,\delta}(C_i))s'}{\rho'}_{\min}^{-1} \delta^{P_{t,\delta}(C_i) t} \ge {\rho'}_{\min}^{-1}  (P_t (C_i) - \epsilon)t + (1- (P_t (C_i) - \epsilon))s'.$$
Thus from the above and the fact that the lower box dimension is stable under closure, it is clear that 
$\underline{\dim}_B (F_i) = \underline{\dim}_B (O_i) 
\ge (P_t (C_i)t + (1- (P_t (C_i)))s'.$ This completes the proof.
\end{proof}

\subsection{Continuity of dimension:}
Here we target to make a note on the continuity of box and Hausdorff dimensions with respect to the well-known Hausdorff metric. In this regard, let us first see that these dimensions are not continuous in general.
\begin{example}
    We choose $A_n= \Big[0, \frac{1}{n}\Big]$ and $A=\{0\}.$ Then one may check that $A_n \to A$ with respect to the Hausdorff metric, but $\dim_B(A_n)=\dim_H(A_n)=1  \nrightarrow \dim_B(A)=\dim_H(A)=0.$
\end{example}
Priyadarshi \cite{Pr3} showed the continuity of the Hausdorff dimension for a class of fractal sets that contains self-similar and self-conformal sets. In the current subsection, we will show the continuity of fractal dimensions for inhomogeneous GD attractors.
\begin{lemma}\label{new17}
    Let $(f_n)$ be a sequence of similarity mappings with similarity ratio $\rho_n$. If $(f_n) $ converges to a function $f$ uniformly then $f$ will be a similarity with similarity ratio $\rho:=\lim_{n \to \infty} \rho_n.$
\end{lemma}
\begin{proof}
    Let us denote the Euclidean norm by $\|.\|_2.$ Using the assumptions taken and continuity of $\|.\|_2$, we have
    \[
    \|f(x) -f(y) \|_2 = \|\lim_{n \to \infty} f_n(x) - \lim_{n \to \infty}f_n(y) \|_2= \lim_{n \to \infty}\| f_n(x)- f_n(y) \|_2 = \lim_{n \to \infty} \rho_n \| x-y\|_2,
    \]
    establishing the claim.
\end{proof}

\begin{theorem}
    Let $\mathcal{GD}_n=(V,E^*,i,t,(X_{v_j},d_{v_j},C_j)_{v_j \in V},(f_e^n)_{e \in E^1})$ be a sequence of GDIFSs consisting of similarity maps $f_e^n$. If each $f_e^n$ converges to $f_e$ and, both $\mathcal{GD}_n$ and $\mathcal{GD}=((f_e)_{e \in E^1})$ with the same vertices and edges satisfy the GDIOSC for all, $n \in \mathbb{N}$ then \begin{itemize}
        \item $\overline{\dim}_B(F_i^n) \to \overline{\dim}_B(F_i)$ provided that $\overline{\dim}_B(C_i)=\overline{\dim}_B(C_j)$ for all $1\le i,j \le N.$
        \item $\dim_H(F_i^n) \to \dim_H(F_i)$ provided that $\dim_H(C_i)=\dim_H(C_j)$ for all $1\le i,j \le N.$
    \end{itemize}
\end{theorem}
\begin{proof}
    Since each $\mathcal{GD}_n$ satisfies GDIOSC and mappings $f_e^n$ are similarity, Corollary \ref{result20} yields $\overline{\dim}_B(F_i^n)=\max\{s_n^*, \overline{\dim}_B(C_i)\}$, where $s_n^*$ is the corresponding Mauldin-Williams graph dimension for $\mathcal{GD}_n.$ Now, Lemma \ref{new17} yields that each $f_e$ is a similarity with similarity ratio $\rho(e):=\lim_{n \to \infty} \rho_n(e)$ because $(f_e^n)$ is a sequence of similarity mappings with similarity ratio $\rho_n(e)$ such that $f_e^n \to f_e$. In view of \cite[Lemma 3.1]{Pr3} and \cite[Lemma 3.2]{Pr3}, we deduce that $s_n^* \to s^*.$ Further, it follows that \[
    \overline{\dim}_B(F_i^n)=\max\{s_n^*,\overline{\dim}_B(C_i)\} \to \max\{s^*,\overline{\dim}_B(C_i)\} =\overline{\dim}_B(F_i).
    \]
    Now, let us turn to the second part. Using Theorem \ref{new615}, we have $$\dim_H (F_i^n) = \max \{\dim_H (F_{\phi_i}^n), \dim_H (C_i)\} \text{ and } \dim_H(F_i) = \max \{\dim_H(F_{\phi_i}), \dim_H(C_i)\}.$$ From the first part of this proof, we get  
    $\dim_H(F_i^n) \to \dim_H(F_i),$ concluding the proof.
\end{proof}
\begin{example}
   In general, if the limit IFS $\mathcal{I}$ does not satisfy the OSC, then the continuity of dimension may not hold. For $n \ge 2,$ choose IFS $\mathcal{I}_n =\bigg\{f_{e_1}^n,f_{e_2}^n: f_{e_1}^n(x)= \frac{x}{n}, ~~f_{e_2}^n(x)= \frac{(n-1)x}{n}+\frac{1}{n}   \bigg\}$ then the homogeneous attractor $F_n,$ corresponding to each $\mathcal{I}_n$ is $[0,1].$ Clearly, $\{\mathcal{I}_n\}$ converges to $ \mathcal{I}= \{f_{e_1},f_{e_2}: f_{e_1}(x)=0 ~~\text{and}~~ f_{e_2}(x)=x\}$ with attractor $F =\{0\} $ but $\dim_H(F_n) \nrightarrow \dim_H(F)$ as $\dim_H(F_n) = 1$ for all $n,$ and $\dim_H(F)=0.$ Here each $\mathcal{I}_n$ for $n \ge 2$ satisfies the OSC but $\mathcal{I}$ does not satisfy the OSC.
\end{example}

\section{Conclusion and future remarks}\label{section 8}
In this article, we have defined the various separation conditions for inhomogeneous GDIFS and studied the dimension of inhomogeneous GDIFS attractors, especially countably unstable dimensions (upper box dimension and lower box dimension). We have shown that inhomogeneous GDIFS also follows a similar dimensional result as in \cite{Fraser1} under inhomogeneous GD open set conditions.
As our given construction effectively generalizes many previous well-known methods such as IFS, graph-directed and inhomogeneous IFS. In the future, we plan to apply this construction in the theory of fractal functions. We also try to estimate dimensions of inhomogeneous GD fractal measures, following the work of Olsen \cite{Olsen1} and Selmi \cite{Selmi1}.

\section*{Statements and Declarations}
\textbf{Data availability:} Data sharing not applicable to this article as no datasets were generated or analysed during the current study.\\
  \textbf{Funding:} The first author thanks IIIT Allahabad (MHRD) for financial support in the form of a Junior Research Fellowship.\\
 \textbf{Conflict of interest:} 
We declare that we do not  have any conflict of interest.\\
 \textbf{Author Contributions:} 
 All authors contributed equally to this manuscript.

    \bibliographystyle{amsplain}

\end{document}